\documentclass[11pt]{article}
\usepackage{latexsym,amssymb,amsmath,amsfonts,amsthm,amsopn,amscd}
\usepackage{graphicx} 
\newtheorem{lem}{Lemma}[section]
\newtheorem{thm}{Theorem}[section]
\newtheorem{cor}{Corollary}[section]

\newcommand{\f}[1]{\mathfrak{#1}}

\newcommand{\up}{\widetilde{Sp}(n, \mb R)}
\newcommand{\mb}{\mathbb}
\newcommand{\commentout}[1]{}

\newcommand{\mc}{\mathcal}
\newcommand{\arr}[1]{\left( \begin{array}{clcr} #1 \end{array} \right)}

\newcommand{\xin}{Sp(n, \mb R)}

\newcommand{\diag}{{\rm diag}}
\begin{document}
\title{Restrictions of the Complementary Series of the Universal Covering of the Symplectic Group }
\author{Hongyu He \\
Department of Mathematics,\\
 Louisiana State University, \\
 Baton Rouge, LA 70803, U.S.A.\\
email: livingstone@alum.mit.edu
}
\date{}
\maketitle
\abstract{In this paper, we study the restrictions of the complementary series representation onto a symplectic subgroup no bigger than half of the size of the original symplectic group.}
\section{Introduction}
Let $\widetilde{Sp}(n,\mb R)$ be the universal covering of $Sp(n,\mb R)$. $\up$ is a central extension of $\xin$:
$$ 1 \rightarrow C  \rightarrow \up \rightarrow \xin \rightarrow 1,$$
where $C \cong \mathbb Z$.
The unitary dual of $C $ is parametrized by a torus $\mathbb T$.   For each $\kappa \in \mathbb T$, denote the corresponding unitary character of $C$ by $\chi^{\kappa}$. We say that a representation $\pi$ of $\up$ is of class $\kappa$ if $\pi|_{C}=\chi^{\kappa}$. Since $C$ commutes with $\widetilde{Sp}(n, \mb R)$, for any irreducible representation $\pi$ of $\widetilde{Sp}(n, \mb R)$, $\pi|_{C}=\chi^{\kappa}$ for some $\kappa$. \\
\\
Denote the projection $\up \rightarrow \xin$ by $p$. For any subgroup $H$ of $\xin$, denote the full inverse image $p^{-1}(H)$ by $\tilde H$. We adopt the notation from ~\cite{sahi}.
Let $P$ be the Siegel parabolic subgroup of $\xin$. One dimensional characters of $\tilde P$ can be parametrized by $(\epsilon, t)$ where $\epsilon \in \mathbb T$ and $t \in \mathbb C$. Let $I(\epsilon, t)$ be the representation of $\up$ induced from the one dimensional character parametrized by $(\epsilon, t)$ of $\tilde P$. If $t$ is purely imaginary, $I(\epsilon, t)$ is unitary and irreducible. If $t$ is real, then $I(\epsilon,t)$ has an invariant Hermitian form. Sahi gives a classification of all irreducible unitarizable $I(\epsilon,t)$. Besides the unitary principal series, there are complementary series $C(\epsilon, t)$ for $t$ in a suitable interval (~\cite{sahi}).\\
\\
Let $(Sp(p, \mb R), Sp(n-p, \mb R))$ be a pair of symplectic groups diagonally embedded in $Sp(n, \mb R)$. Suppose that $p \leq n-p$. Let $U(n)$ be a maximal compact subgroup such that $Sp(n-p, \mb R) \cap U(n)$ is a maximal compact subgroup of $Sp(n-p, \mb R)$. Denote $Sp(n-p, \mb R) \cap U(n)$ by $U(n-p)$.
\begin{thm} Suppose that $p \leq n-p$ and $C(\epsilon,t)$ is unitary. Then
$$C(\epsilon, t)|_{\tilde U(n-q) {\widetilde Sp}(p, \mb R)} \cong I(\epsilon, 0)|_{\tilde U(n-q) {\widetilde Sp}(p, \mb R)} \cong I(\epsilon, i \lambda)|_{\tilde U(n-q) {\widetilde Sp}(p, \mb R)} \qquad (\lambda \in \mb R).$$
\end{thm}
$p=[\frac{n}{2}]$ is the best possible value for such a statement. In particular, for $\widetilde{Sp}(2m+1, \mb R)$
$$ I(\epsilon, 0)|_{\widetilde{Sp}(m+1, \mb R)} \ncong C(\epsilon, t)|_{\widetilde{Sp}(m+1, \mb R)}.$$
To see this, let $L^2(\up)_{\kappa}$ be the set of functions with
$$f(zg)=\mu^{\kappa}(z) f(g) \qquad (z \in C, g \in \up);$$
$$\| f \|^2=\int_{\xin} |f(g)|^2 d [g] \qquad (g \in \up, [g] \in \xin).$$
We say that a representation of class $\kappa$ is tempered if it is weakly contained in $L^2(\up)_{\kappa}$. By studying the leading exponents of $I(\epsilon,0)$ and $C(\epsilon, t)$, it can be shown that
$I(\epsilon, 0)|_{\widetilde{Sp}(m+1, \mb R)}$ is \lq\lq tempered \!\rq\rq and $C(\epsilon, t)|_{\widetilde{Sp}(m+1, \mb R)}$ is not \lq\lq tempered \!\rq\rq. Therefore 
$$ I(\epsilon, 0)|_{\widetilde{Sp}(m+1, \mb R)} \ncong C(\epsilon, t)|_{\widetilde{Sp}(m+1, \mb R)}.$$
The author would like to thank G. Olafsson and J. Lawson for some very help discussions.

\section{A Lemma on Friedrichs Extension}
Let $S$ be a semibounded densely defined symmetric operator on a Hilbert space $H$. Suppose that $(Su, u) > 0$ for every nonzero $u \in \mathcal D(S)$. We call $S$ positive. For $u,v \in \mathcal D(S)$,
define
$$(u,v)_S=(u, S v),$$
$$\|u\|_S=(u, Su).$$
Let $H_S$ be the completion of  $\mc D(S)$ under the norm $\|\ \|_S$.\\
\\
The operator $S+I$ has a unique Friedrichs extension $(S+I)_0$ such that $\mc D((S+I)_0) \subseteq H_{S+I}$ and $(u, v)_{S+I}=(u, (S+I)_0 v)$ for all $u \in H_{S+I}$ and $v \in \mathcal D((S+I)_0)$ (see Theorem in Page 335 ~\cite{rs}).  Here $H_{S+I} \subseteq H$ and $(S+I)_0$ is self-adjoint. Now consider $(S+I)_0-I$. It is an self-adjoint extension of $S$. It is nonnegative. By the spectral decomposition and functional calculus, $(S+I)_0-I$ has unique square root $T$ (See 127. 128. ~\cite{rs}).
\begin{lem} Let $S$ be a positive densely defined symmetric operator. Then the square root of $(S+I)_0-I$ extends to an isometry from $H_S$ into $H$. 
\end{lem}
Proof: Clearly, the spectrum of $T$ is contained in the nonnegative part of the real line. By spectral decomposition $\mc D((S+I)_0-I)=\mc D((S+I)_0) \subseteq \mc D(T)$ and $T T=(S+I)_0-I$.  In addition 
for any $u,v \in \mc D(S) \subseteq D((S+I)_0)$, 
$$(Tu, T v)=(u, TT v)=(u, (S+I)_0 v-v)=(u, S v)=(u,v)_S.$$
So $T$ is an isometry from $\mc D(S)$ into $H$. Since $\mc D(S)$ is dense in $H_S$, $T$ extends to an isometry from $H_S$ into $H$. $\Box$ \\
\\
We denote the isometry by $I_S$. It is canonical. 
\section{Complementary Series of $\up$}
Fix the Lie algebra $\f{sp}(n, \mb R)$:
$$\{ \arr{ X & Y \\ Z & -X^t} \mid Y^t=Y, Z^t=Z \}$$
and the Siegel parabolic algebra $\f p$:
$$\{ \arr{ X & Y \\ 0 & -X^t} \mid Y^t=Y \}.$$
Fix the Levi decomposition $\f p=\f l \oplus \f n$ with  
$$\f l= \{ \arr{ X & 0 \\ 0 & -X^t } \mid X \in \f{gl}(n,\mb R) \}$$
and 
$$\f n =\{ \arr{ 0 & Y \\ 0 & 0} \mid Y^t=Y \}.$$
Fix a Cartan  subalgebra 
$$\f a= \{ \diag(H_1,H_2, \ldots, H_n, -H_1, -H_2 \ldots, -H_n) \mid H_i \in \mathbb R \}.$$
Let $\xin$ be the symplectic group and $P$ be the Siegel parabolic subgroup. Let $LN$ be the Levi decomposition and $A$ be the analytic group generated by the Lie algebra $\f a$. Clearly, $L \cong GL(n, \mb R)$ and $L \cap U(n) \cong O(n)$. On the covering group, we have $\tilde L \cap \tilde U(n)=\tilde O(n)$. Recall that
$$\tilde U(n)= \{ (x, g) \mid g \in U(n), \exp 2 \pi i x = \det g, x \in \mathbb R \}.$$
Therefore
$$\tilde O(n)=\{ (x, g) \mid g \in O(n), \exp 2 \pi i x= \det g, x \in \mathbb R \}.$$
Notice that $ x \in \frac{1}{2} \mathbb Z$ since $\det g= \pm 1$.
We have the following exact sequence
$$ 1 \rightarrow SO(n) \rightarrow \tilde O(n) \rightarrow \frac{1}{2} \mathbb Z \rightarrow 1.$$
Consequently, we have
$$ 1 \rightarrow GL_0(n, \mb R) \rightarrow \tilde L \rightarrow \frac{1}{2} \mathbb Z \rightarrow 1.$$
In fact, $$\tilde L=\{ (x, g) \mid g \in L, \exp 2 \pi i x= \frac{\det g}{|\det g|}, x \in \mathbb R \}.$$
The one dimensional unitary characters of $\frac{1}{2} \mathbb Z$ are parametrized by the one dimensional torus $T$. Identify $T$ with $[0, 1)$. Let $\mu^{\epsilon}$ be the character of $\frac{1}{2} \mathbb Z$ corresponding to $\epsilon \in [0,1)$
Now each character $\mu^{\epsilon}$ yields a character of $\tilde L$, which in turn, yields a character of $\tilde P$. For simplicity, we retain $\mu^{\epsilon}$ to denote the character on $\tilde L$ and $\tilde P$. Let $\nu$ be the $\det$-character on ${\tilde L}_0$, i.e.,
$$ \nu(x,g)= | \det g | \qquad (x,g) \in \tilde L.$$
Let
$$I(\epsilon, t)=Ind_{\tilde P}^{\up} \mu^{\epsilon} \otimes \nu^{t}$$
be the normalized induced representation with $\epsilon \in [0,1)$ and $t \in \mathbb C$. 
This is Sahi's notation applied to the universal covering of the symplectic group ~\cite{sahi}.  $I(\epsilon,t)$ is a degenerate principal series representation. 
Clearly, $I(\epsilon, t)$ is unitary when $t$ is purely imaginary.\\
\\
When $t$ is real and $I(\epsilon, t)$ is unitarizable, the unitary representation, often denoted by $C(\epsilon, t)$, is called a complementary series representation. 
Various complementary series of  $\xin$ and its metaplectic covering was determined explicitly or implicitly by  Kudla-Rallis, {\O}rsted-Zhang, Brason-Olafsson-{\O}rsted and others. See
 ~\cite{kr}, ~\cite{boo}, ~\cite{oz} and the references therein. The complete classification of the complementary series of the universal covering is due to Sahi (~\cite{sahi}).
\begin{thm}[Thm A, ~\cite{sahi}]
Suppose that $t$ is real.
For $n$ even, $I(\epsilon,t)$ is irreducible and unitarizable if and only if 
$0< |t| <|\frac{1}{2}-|2 \epsilon-1||$. For $n$ odd and $n >1$, $I(\epsilon, t)$ is irreducible and unitarizable if and only if $0< |t| < \frac{1}{2}-|\frac{1}{2}-|2 \epsilon-1||$.
\end{thm}
\begin{figure}[htbp]
	\centering
		\caption{Complementary Parameters $(E,t)$}
	\label{fig:Diagram1}
\end{figure}
 One can easily check that the complementary series exist if $\epsilon \neq 0, \frac{1}{2}$ for $n$ odd and $n>1$ ; if $\epsilon \neq \frac{1}{4}, \frac{3}{4}$ for $n$ even. It is interesting to note that complementary series always exist unless $I(\epsilon, t)$ descends into a representation of the metaplectic group. For the metaplectic group ${Mp}(2n+1, \mb R)$, there are two complementary series $I(\frac{1}{4}, t) ( 0 < t < \frac{1}{2})$ and $I(\frac{3}{4}, t) ( 0< t< \frac{1}{2})$. For the metaplectic group $Mp(2n, \mb R)$, there are two complementary series
$I(0, t) ( 0 < t < \frac{1}{2})$ and $I(\frac{1}{2}, t) ( 0< t< \frac{1}{2})$.
These four complementary series are the \lq\lq longest \rq\rq.\\
\\
For $n=1$, the situation is quite different. The difference was  pointed out in ~\cite{kr}. For example, there are Bargmann's complementary series representation for $I(0, t) (t \in (0, \frac{1}{2}))$. The classification of the complementary series of  $\widetilde{Sp}(1, \mb R)$ can be found in ~\cite{bar}, ~\cite{puk}, ~\cite{howe}.\\
\\
Since our restriction theorem only makes sense for $n \geq 2$. We will assume $n \geq 2$ from now on. The parameters for the complementary series of $\widetilde{Sp}(n, \mb R)$ are illustrated in fig. \! ~\ref{fig:Diagram1}.

\section{The generalized compact model and The Intertwining Operator}
Recall that $$I^{\infty}(\epsilon, t)= \{ f \in C^{\infty}(\up) \mid f(g l n)= (\mu^{\epsilon} \otimes \nu^{t+\rho})(l^{-1}) f (g) ( l \in \tilde L, n \in N) \}$$
where $\rho=\frac{n+1}{2}$.
Let $X= \up/\tilde P$.  Then $I^{\infty}(\epsilon, t)$ consists of smooth sections of
the homogeneous vector bundle $\mc L_{\epsilon, t}$
$$ \up \times_{\tilde P} \mathbb C_{\mu^{\epsilon} \otimes \nu^{t+\rho}} \rightarrow X.$$
Since $X \cong \xin/P \cong \tilde U(n)/\tilde O(n)$, $\tilde U(n)$ acts transitively on $X$. $f \in I^{\infty}(\epsilon, t)$ is uniquely determined by $f|_{\tilde U(n)}$ and vice versa. Moreover, the homogeneous vector bundle
$\mc L_{\epsilon, t}$ can be identified with $\mc K_{\epsilon, t}$
$$ \tilde U(n) \times_{\tilde O(n)} C_{\mu^{\epsilon} \otimes \nu^{t+\rho}}|_{\tilde O(n)} \rightarrow X$$
naturally. Notice that the homogeneuous vector bundle $\mc K_{\epsilon, t}$ does not depend on the parameter $t$. We denote it by $\mc K_{\epsilon}$. The representation $I(\epsilon, t)$ can then be modeled on smooth sections of $\mc K_{\epsilon}$. This model will be called
the generalized compact model.\\
\\
The generalized compact model provides much convenience. First, it equips the smooth sections of $\mc K_{\epsilon, t}$ with a pre-Hilbert structure
$$(f_1, f_2)_X= \int_{[k] \in X}  f_1(k) \overline{f_2(k)} d [k],$$
where $k \in \tilde U(n)$ and $[k] \in X$.
It is easy to verify that $f_1(k) \overline{ f_2(k)}$ is a function of $[k]$ and it does not depend on any particular choice of $k$.
Notice that our situation is different from the compact model since $\tilde U(n)$ is not compact. We denote the completion of $I^{\infty}$ with respect to $(,)_X$ by $I_X(\epsilon, t)$.  Secondly, the action of $\tilde U(n)$ on
$\mc K_{\epsilon}$ induces an orthogonal decomposition of $I_X(\epsilon,t)$:
$$I_X(\epsilon, t)= \hat{\oplus}_{\alpha \in 2 \mathbb Z^n} V( \alpha + \epsilon(2,2,\ldots, 2)),$$
where $V(\alpha + \epsilon(2,2,\ldots, 2))$ is an irreducible representation of $\tilde U(n)$ with highest weight $\alpha+ \epsilon(2,2,\ldots, 2)$. Let
$$V(\epsilon, t)= {\oplus}_{\alpha \in 2 \mathbb Z^n} V( \alpha + \epsilon(2,2,\ldots, 2)).$$
$V(\epsilon, t)$ possesses an action of the Lie algebra $\f{sp}(n, \mathbb R)$. It is called the Harish-Chandra module of $I(\epsilon, t)$. Clearly, $V(\epsilon, t) \subset I^{\infty}(\epsilon,t) \subset I_X(\epsilon, t)$. \\
\\
For each $t$, there is an $\up$-invariant sesquilinear pairing of $I(\epsilon, t)$ and $I(\epsilon, -\overline{t})$, namely,
$$(f_1, f_2)=\int_{X} f_1(k) \overline{f_2(k)} d [k],$$
where $f_1 \in I(\epsilon, t)$ and $f_2 \in I(\epsilon, -\overline{t})$. If $t$ is purely imaginary,
we obtain a $\up$-invariant Hermitian form which is exactly $(,)_X$. Since $(,)_X$ is positive definite, $I(\epsilon, t)$ is unitary.\\
\\
For each real $t$, the form $(,)$ gives an $\up$-invariant sesquilinear pairing of $I(\epsilon, t)$ and $I(\epsilon, -t)$. 
There is an intertwining operator 
$$A(\epsilon, t): V(\epsilon, t) \rightarrow V(\epsilon, -t)$$
which preserves the action of $\f{sp}(n, \mathbb R)$ (see for example ~\cite{boo}). Define a Hermitian structure
$(,)_{\epsilon, t}$ on $V(\epsilon, t)$ by 
$$ (u, v)_{\epsilon, t}=(A(\epsilon, t) u, v), \qquad (u, v \in V(\epsilon, t)).$$
Clearly, $(,)_{\epsilon, t}$ is $\f{sp}(n, \mb R)$-invariant. So $A(\epsilon, t)$ induces an invariant Hermitian form on
$V(\epsilon, t)$. \\
\\
Now $A(\epsilon, t)$ can also be realized as an unbounded operator on $I_{X}(\epsilon, t)$ as follows. For each $f \in V(\epsilon, t)$, define $A_X(\epsilon, t) f$ to be the unique section of $\mc L_{\epsilon, t}$ such that
$$(A_X(\epsilon, t) f) |_{\tilde U(n)}= (A(\epsilon, t)f)|_{\tilde U(n)}.$$
Notice that $A_X(\epsilon,t)f \in I(\epsilon, t)$  and $A(\epsilon, t) f \in I(\epsilon,-t)$. They differ by a multiplier. \\
\\
Now $A_X(\epsilon, t)$ is an unbounded operator on the Hilbert space $I_{X}(\epsilon, t)$.
 The following fact is well-known in many different forms. I state it in a way that is convenient for later use.
\begin{lem} Let $t \in \mathbb R$.
$I(\epsilon,t)$ is unitarizable if an only if $A_X(\epsilon, t)$ extends to a self-adjoint operator on $I_X(\epsilon, t)$ with spectrum on the nonnegative part of the real axis.
\end{lem}
The spectrum of $A_{X}(\epsilon, t)$ was computed in ~\cite{boo} and ~\cite{oz} explicitly for special cases and in ~\cite{sahi} implicitly. In particular, $A_X(\epsilon, t)$ restricted onto each $\tilde U(n)$-type is a scalar multiplication and the scalar is bounded by a polynomial on the highest weight. We obtain
\begin{lem}
$A_X(\epsilon, t)$ extends to an unbounded operator from $I^{\infty}(\epsilon,t)$ to $I^{\infty}(\epsilon, t)$.
\end{lem}
This lemma follows from a standard argument that the norm of each $\tilde U(n)$-component in the Peter-Weyl expansion of any smooth section of $\mc K_{\epsilon}$ decays rapidly with respect to the highest weight.
\section{The Mixed Model }
Suppose that $p+q=n$ and $p \leq q$. Fix a subgroup $Sp(p, \mb R) \times Sp(q,  \mb R)$ in $\xin$.  Then we have a subgroup
$\widetilde{Sp}(p, \mb R) \widetilde{Sp}(q, \mb R)$ in $\up$. Notice that $\widetilde{Sp}(p, \mb R) \cap \widetilde{Sp}(q, \mb R) \cong \mathbb Z$. So $\widetilde{Sp}(p, \mb R) \widetilde{Sp}(q, \mb R)$ is not a direct product, but rather the product of the two groups as sets. Let $U(q)=Sp(q, \mb R) \cap U(n)$. 
\begin{thm}[Main Theorem]
Suppose that $p+q=n$ and $p \leq q$.
Given a complementary series representation $C(\epsilon, t)$,
$$C(\epsilon, t)|_{ \widetilde{Sp}(p, \mb R) \tilde U(q)} \cong I(\epsilon, 0)|_{ \widetilde{Sp}(p, \mb R) \tilde U(q)}.$$
In other words, there is an isometry between $C(\epsilon, t)$ and $I(\epsilon, 0)$ that intertwines the actions of $\tilde U(q)$ and of
$\widetilde{Sp}(p, \mb R)$.
\end{thm}
We begin by recall a result concerning the action of $Sp(p, \mb R) \times Sp(q, \mb R)$ on $X$ (~\cite{he06}).
\begin{lem}~\label{geometry} $Sp(p, \mb R) \times Sp(q, \mb R)$ acts on $X$ with a unique open dense orbit $X_0$. Let $p \leq q$. Let $P_{p, 2q-2p}$ be a maximal parabolic subgroup of $Sp(q,\mb R)$ preserving an $q-p$ dimensional isotropic subspace.  Let $GL_{q-p} Sp(p, \mb R) N_{p, 2q-2p}$ be the Langlands decomposition of $P_{p, 2q-2p}$. 
Let $$H= \{ ( u, m_2 \ {}^tu^{-1} n_2) \mid m_2 \in GL(q-p, \mb R), u \in Sp(p, \mb R), n_2 \in N_{p, 2q-2p} \}.$$
Then $X_0 \cong Sp(p, \mb R) \times Sp(q, \mb R)/H$.
\end{lem}
Notice that $H \cong P_{p, 2q-2p}$. But the $Sp(p, \mb R)$ factor in $P$ is diagonally embedded in $Sp(p, \mb R) \times Sp(q, \mb R)$.
In particular, there is a principal fibration
$$ Sp(p, \mb R) \rightarrow X_0 \rightarrow Sp(q, \mb R)/P_{p, 2q-2p} \cong U(q)/O(q-p)U(p).$$
Here $O(q-p)U(p)=U(q) \cap P_{p, 2q-2p} \subset Sp(q, \mb R)$.
Let $M=Sp(p, \mb R) U(q)$ and $\tilde M= \widetilde{Sp}(p, \mb R) \tilde U(q)$. Then
$\mc L_{\epsilon, t}$ restricted onto $X_0$ becomes
$$\tilde M \times_{\widetilde{ O(q-p) U(p)} } \mathbb C_{\mu^{\epsilon}} \rightarrow \tilde M/\widetilde{ O(q-p) U(p)}  \cong X_0.$$
Let $I_{c, X_0}^{\infty}(\epsilon, t)$ be the set of smooth section of $\mc L_{\epsilon, t}$ that are compactly supported on $X_0$. Clearly $$I_{c, X_0}^{\infty}(\epsilon, t) \subset I^{\infty}(\epsilon, t).$$
Consider the restriction of $(,)_X$ onto $I_{c, X_0}^{\infty}(\epsilon, t)$. We are interested in expressing $(,)_X$ as an integral on $\tilde M/\widetilde{ O(q-p) U(p)} $. This boils down to a change of variables from $\tilde U(n)/\tilde O(n)$ to $\tilde M/\widetilde{ O(q-p) U(p)} $.\\
\\
Let $ d g_1$ be the invariant measure on $\widetilde{Sp}(p, \mb R)$ and $d [k_2]$ be the invariant measure on $\tilde U(q)/ \widetilde{ O(q-p) U(p)} $.
Every element in $\up$ has a $\tilde U(n) P_0$ decomposition where $P_0$ is the identity component of $\tilde P$. For each $ g \in \up$, write
$g= \tilde u(g) p(g)$. For each $g_1 k_2 \in \tilde M$, write $g_1 k_2= \tilde u(g_1 k_2) p(g_1 k_2)$. Then $\tilde u$ defines a map from $\tilde M$ to $\tilde U(n)$. $\tilde u$ induces a map from $\tilde M/\widetilde{ O(q-p) U(p)} $ to $\tilde U(n)/\tilde O(n)$ which will be denoted by $j$. Clearly, $j$ is an injection from $\tilde M/\widetilde{ O(q-p) U(p)} $ onto $X_0$. 
Change the variable on $X$ from $\tilde M/ \widetilde{ O(q-p) U(p)} $ to $ \tilde U(n)/ \tilde O(n)$. Let $J([g_1], [k_2])$ be the Jacobian:
$$\frac{ d j([g_1], [k_2])}{ d [g_1] d [k_2]}.$$
 We have
\begin{lem}~\label{2}
Let $$\Delta_{\epsilon,t}(g_1, k_2)=\nu(p(g_1 k_2))^{t+\overline{t}+ 2 \rho}J([g_1], [k_2]).$$ Then for every $f_1, f_2 \in I^{\infty}(\epsilon, t)$ we have
$$(f_1, f_2)_X= \int_{\tilde M/\widetilde{ O(q-p) U(p)} } f_1(g_1 k_2) \overline{f_2(g_1 k_2)} \Delta_{\epsilon, t}(g_1, k_2) d [g_1] d [k_2]$$
where $g_1 \in \widetilde{Sp}(p, \mb R)$, $k_2 \in \tilde U(q)$, $[g_1] \in Sp(p,\mb R)$ and $[k_2] \in \tilde U(q)/\widetilde{ O(q-p) U(p)} $. 
Furthermore,  $\Delta_{\epsilon, t}(g_1, k_2)$ is a nonnegative right $\widetilde{ O(q-p) U(p)} $-invariant function  on 
$\widetilde{Sp}(p, \mb R) \times \tilde U(q)$.
\end{lem}
Proof: We compute
\begin{equation}
\begin{split}
 & \int_{\tilde M/\widetilde{ O(q-p) U(p)} } f_1(g_1 k_2) \overline{f_2(g_1 k_2)} \Delta_{\epsilon, t}(g_1, k_2) d [g_1] d [k_2] \\
 = & \int_{\tilde M/\widetilde{ O(q-p) U(p)} } f_1(\tilde u(g_1 k_2)) \overline{f_2(\tilde u(g_1 k_2))} \nu(p(g_1 k_2))^{-t-\overline{t}- 2 \rho} \Delta_{\epsilon, t}(g_1, k_2) d [g_1] d [k_2] \\
 = & \int_{\tilde M/\widetilde{ O(q-p) U(p)} } f_1(\tilde u(g_1 k_2)) \overline{f_2(\tilde u(g_1 k_2))} \nu(p(g_1 k_2))^{-t-\overline{t}- 2 \rho} \Delta_{\epsilon, t}(g_1, k_2) J^{-1}(g_1, k_2) d j([g_1][ k_2]) \\
 = & \int_{X_0} f_1(\tilde u) \overline{f_2(\tilde u)} d [\tilde u]=(f_1, f_2)_X.
 \end{split}
 \end{equation}
Since $\nu(p(g_1 k_2)$ and $J([g_1], [k_2])$ remains the same when we multiply $k_2$ on the right by $\widetilde{ O(q-p) U(p)} $, $\Delta_{\epsilon, t}(g_1, k_2)$ is a nonnegative right $\widetilde{ O(q-p) U(p)} $-invariant function. $\Box$ \\
\\
For each $f_1, f_2 \in I^{\infty}_{c, M}(\epsilon, t)$, define 
$$(f_1 , f_2)_{M, t}= \int_{M} f_1(g_1 k_2) \overline{f_2(g_1 k_2)} \Delta_{\epsilon, t}(g_1 k_2) d [g_1] d [k_2],$$
$$(f_1, f_2)_{M}= \int_{M} f_1(g_1, k_2) \overline{f_2(g_1 k_2)} d [g_1] d [k_2].$$
Let $I_M(\epsilon, t)$ be the completion of $I^{\infty}_{c, M}(\epsilon, t)$ under $(,)_{M, t}$. Clearly $I_M(\epsilon, t) \cong I_X(\epsilon, t)$ as Hilbert representations of $\up$. $I_M(\epsilon, t)$ is called the mixed model.
 \section{Mixed Model for Unitary Principal Series}
\begin{lem}~\label{1} If $t$ is purely imaginary, then $\Delta_{\epsilon, t}(g_1, k_2 )$ is a constant and $(,)_{M,t}$ is a constant multiple of $(,)_{M}$.
\end{lem}
Proof: Let $t \in i \mb R$. Let $f_1, f_2 \in I^{\infty}(\epsilon, t)$ and $h \in \widetilde{Sp}(p, \mb R)$. We have
\begin{equation}
\begin{split}
& (I(\epsilon, t)(h) f_1, I(\epsilon, t)(h) f_2)_X \\
= & (I(\epsilon, t)(h) f_1, I(\epsilon, t)(h) f_2)_{M,t} \\
= & \int_{\tilde M/\widetilde{ O(q-p) U(p)} } f_1(h^{-1}g_1 k_2) \overline{f_2(h^{-1} g_1 k_2)} \Delta_{\epsilon t}(g_1, k_2) d [g_1] d [k_2]\\
= & \int_{\tilde M/\widetilde{ O(q-p) U(p)} } f_1(g_1 k_2) \overline{f_2(g_1,k_2)} \Delta_{\epsilon, t}(h g_1, k_2) d [g_1] d [k_2] \\
\end{split}
\end{equation}
Since $I(\epsilon, t)$ is unitary, $$(I(\epsilon, t)(h) f_1, I(\epsilon, t)(h) f_2)_X=
(f_1, f_2)_X.$$ 
We have 
$$\int_{\tilde M/\widetilde{ O(q-p) U(p)} } f_1(g_1 k_2) \overline{f_2(g_1 k_2)} \Delta_{\epsilon, t}(h g_1, k_2) d [g_1] d [k_2]=
\int_{\tilde M/\widetilde{ O(q-p) U(p)} } f_1(g_1 k_2) \overline{f_2(g_1 k_2)} \Delta_{\epsilon, t}(g_1, k_2) d [g_1] d [k_2].$$
It follows that
$\Delta_{\epsilon, t}(h g_1, k_2)=\Delta_{\epsilon, t}(g_1, k_2).$
 for any $h \in \widetilde{Sp}(p, \mb R)$. Similarly, we obtain 
$\Delta_{\epsilon, t}(g_1, k k_2  )=\Delta(g_1, k_2)$
for any $k \in \tilde U(q)$. Hence, $\Delta_{\epsilon, t}(g_1, k_2)$ is a constant for purely imaginary $t$. $\Box$
\begin{cor}~\label{jacobianandnu} ${J([g_1], [k_2])}= c \nu(p(g_1 k_2))^{ -2 \rho}$ and 
$\Delta_{\epsilon,t}(g_1, k_2)=c \nu(p(g_1 k_2))^{t+\overline{t}}$. Furthermore, 
\begin{equation}~\label{mix}
I_M(\epsilon, t) \cong L^2(\tilde M \times_{\widetilde{ O(q-p) U(p)} } \mathbb C_{\mu^{\epsilon}}, \nu(p(g_1 k_2))^{t+\overline{t}} d [g_1] d [k_2]).
\end{equation}
\end{cor}
From now on, identify $L^2(\tilde M \times_{\widetilde{ O(q-p) U(p)} } \mathbb C_{\mu^{\epsilon}}, \nu(p(g_1 k_2))^{t+\overline{t}} d [g_1] d [k_2])$ with $I_M(\epsilon, t)$. 
\begin{cor}~\label{6.2} $\nu(p(g_1 k_2))^{-\rho} \in L^2(\tilde M/\widetilde{ O(q-p) U(p)} )$ and $\nu(p(g_1 k_2))^{-1}$ is a bounded positive function.
\end{cor}
Proof: Since $X$ is compact, $J([g_1], [k_2]) \in L^1(\tilde M/\widetilde{ O(q-p) U(p)} )$. Notice that $\nu(p(g_1 k_2))^{-2 \rho}= c {J([g_1], [k_2])}$. It follows that $\nu(p(g_1 k_2))^{-\rho} \in L^2(\tilde M/\widetilde{ O(q-p) U(p)} )$.  Since ${J([g_1], [k_2])}$ is continuous and $[k_2]$ is a compact manifold, it suffice to show that for each $k_2$, $\{J([g_1], [k_2])\}_{g_1 \in Sp(p, \mb R)}$ is bounded. This is true because $j|_{k_2 Sp(p, \mb R)}$ is an analytic compactification of $Sp(p, \mb R)$ and the Jacobian can be easily computed just like the way it is done in ~\cite{he060} (see Theorem 2.3). In particular, it is always positive. By Cor.
~\ref{jacobianandnu}, $\nu(p(g_1 k_2))^{-1}$ is bounded and positive.  $\Box$\\
\\
If $f \in I_{M}(\epsilon, t_1)$ and $h > 0$, we have
$ \|f \|_{M, t_1} \geq C  \|f \|_{M, t_1-h}$. So $I_{M}(\epsilon, t_1) \subset I_{M}(\epsilon, t_1-h)$. 
\begin{cor} Suppose that $h >0$. Then $I_{M}(\epsilon, t_1) \subset I_{M}(\epsilon, t_1-h)$ under the Equation \!~\ref{mix}.
\end{cor}
Notice that in the mixed model, the actions of $\widetilde{Sp(p,\mb R)}$ and $\tilde U(q)$ does not depend on the parameter $t$. We obtain
\begin{thm}~\label{main0}
Let $t$ be purely imaginary. $I(\epsilon, t)$ can all be modeled on  $$L^2(\tilde M \times_{\widetilde{ O(q-p) U(p)} } \mathbb C_{\mu^{\epsilon}}, d [g_1] d [k_2]).$$ In particular,
$I_M(\epsilon, t)|_{\widetilde{Sp}(p, \mb R) \tilde U(q)} \cong I_M(\epsilon, 0)|_{\widetilde{Sp}(p, \mb R) \tilde U(q)}$ and the identity operator intertwines $I_M(\epsilon,0)|_{\widetilde{Sp}(p, \mb R) \tilde U(q)}$ with $I_M(\epsilon, t)|_{\widetilde{Sp}(p, \mb R) \tilde U(q)}$.
\end{thm}
For $t$ a nonzero real number, $\Delta_{\epsilon, t}(g,k)$ is not a constant. So Theorem ~\ref{main0} does not hold for real nonzero $t$. \\
\\
\section{\lq\lq Square Root \!\rq\rq of the Intertwining Operator}

Suppose from now on $t \in \mathbb R$. For $ f \in I^{\infty}(\epsilon, t)|_{\tilde M}$, define a function on $\tilde M$,
$$(A_M(\epsilon, t) f)(g_1 k_2)= A(\epsilon, t) f (g_1 k_2) \qquad (g_1 \in \widetilde{Sp}(p, \mathbb R), k_2 \in \tilde U(q)).$$
So $A_M(\epsilon, t)$ is the \lq\lq restriction \rq\rq of $A(\epsilon, t)$ onto $\tilde M$. $A_M(\epsilon,t)$ is not yet an unbounded operator on $I_M(\epsilon, t)$. In fact, for $t >0$, $A_M(\epsilon, t)$ does not behave well. In this case, it is not clear whether
$A_M(\epsilon,t)$ can be realized as an unbounded operator on $I_M(\epsilon, t)$. 
$A_M(\epsilon,t) f$ differs from $A_X(\epsilon, t) f$ by a multiplier.
\begin{lem}~\label{7.1} For $t \in \mb R$ and $f \in I^{\infty}(\epsilon, t)$,
$$(A_M(\epsilon, t) f|_{\tilde M})(g_1 k_2)=(A_X(\epsilon, t) f)(g_1 k_2) \nu(p(g_1k_2))^{2t}=(A_X(\epsilon, t) f)(g_1 k_2) \Delta_{\epsilon, t}(g_1, k_2).$$
\end{lem}
This Lemma is due to the fact that $A_X(\epsilon,t)f \in I(\epsilon,t)$ but $A(\epsilon,t) f \in I(\epsilon, -t)$.\\
\\
Let $f \in I^{\infty}(\epsilon, t)$. In terms of the mixed model, the invariant Hermitian form $(,)_{\epsilon, t}$ can be written as follows:
$$(f,f)_{\epsilon, t}=(A_X(\epsilon,t) f, f)_X=\int_{\tilde M/{\widetilde O(q-p)U(p)}}
A_M(\epsilon, t) \, f|_{\tilde M} \ \overline{f}|_{\tilde M} d [g_1] d [k_2].$$
We obtain
\begin{lem}~\label{7.2} For $f_1, f_2 \in I^{\infty}(\epsilon, t)$, $(f_1,f_2)_{\epsilon, t}=(A_M(\epsilon, t) f_1|_{\tilde M}, f_2|_{\tilde M})_{M}$.
\end{lem}
\begin{thm} If $t <0$ and $C(\epsilon,t)$ is a complementary series representation, then $A_M(\epsilon, t)$ is positive and densely defined symmetric operator. Its self-adjoint-extension $(A_M(\epsilon, t)+I)_0-I$ has a unique square root which extends to an isometry from $C(\epsilon, t)$ onto 
$$L^2(\tilde M \times_{\widetilde{ O(q-p) U(p)} } \mathbb C_{\mu^{\epsilon}}, d [g_1] d [k_2]) .$$ 
\end{thm}
Proof: Let $t <0$. Put 
$$\mc H= L^2(\tilde M \times_{\widetilde{ O(q-p) U(p)} } \mathbb C_{\mu^{\epsilon}}, d [g_1] d [k_2]) .$$ 
Let $f \in I^{\infty}(\epsilon, t)$. Then 
$$A_M(\epsilon, t ) (f|_{\tilde M})=\nu(p(g_1 k_2))^{2t} A_X(\epsilon, t) f(g_1 k_2).$$
By Lemma ~\ref{7.1}, Cor. ~\ref{6.2} and Lemma ~\ref{2}, we have
\begin{equation}
\begin{split}
 &  \int_{\tilde M/{\widetilde O(q-p)U(p)}}
A_M(\epsilon, t) (f|_{\tilde M}) \overline{A_M(\epsilon,t)(f |_{\tilde M})} d [g_1] d [k_2] \\
= &\int_{\tilde M/{\widetilde O(q-p)U(p)}}
\nu(p(g_1 k_2))^{2t} |(A_X(\epsilon, t) f) (g_1 k_2)|^2 \Delta_{\epsilon, t}(g_1, k_2) d [g_1] d [k_2] \\
\leq  & C \int_{\tilde M/{\widetilde O(q-p)U(p)}}
 |A_X(\epsilon, t) f (g_1 k_2)|^2 \Delta_{\epsilon, t}(g_1, k_2) d [g_1] d [k_2] \\
= & C ( A_X(\epsilon, t) f, A_X(\epsilon, t) f)_X < \infty.
\end{split}
\end{equation}
 Therefore, $A_M(\epsilon, t)(f|_{\tilde M}) \in \mc H$. 
Let $\mc D= I^{\infty}(\epsilon, t)|_{\tilde M}$. Clearly, $\mc D$ is dense in $ \mc H$. So $A_M(\epsilon,t)$ is a densely defined unbounded operator. It is positive and symmetric by Lemma ~\ref{7.2}. \\
\\
Now $(f, g)_{\epsilon, t}=(A_M(\epsilon, t)f|_{\tilde M}, g|_{\tilde M})_{M}$ for any $f, g \in I^{\infty}(\epsilon, t)$. So
$C(\epsilon, t) = \mc H_{A_M(\epsilon,t)}$.
By Lemma 2.1, there exists an isometry $I_{A_M(\epsilon, t)}$ mapping from $C(\epsilon, t)$ into $\mc H$. \\
\\
Suppose that
$I_{A_M(\epsilon,t)}$ is not onto. Let $f \in \mc H$ such that for any $u \in \mc D(((A_M(\epsilon,t)+I)_0-I)^{\frac{1}{2}})$,
$$(f, ((A_M(\epsilon,t)+I)_0-I)^{\frac{1}{2}} u)_M=0.$$
Notice that $$I^{\infty}(\epsilon, t)|_{\tilde M} \subset \mc D((A_M(\epsilon,t)+I)_0-I) \subset \mc D(((A_M(\epsilon,t)+I)_0-I)^{\frac{1}{2}}),$$
and
$$(A_M(\epsilon, t)+I)_0-I)^{\frac{1}{2}} (A_M(\epsilon, t)+I)_0-I)^{\frac{1}{2}}=(A_M(\epsilon,t)+I)_0-I.$$
In particular, $$ ((A_M(\epsilon,t)+I)_0-I)^{\frac{1}{2}} I^{\infty}(\epsilon,t)|_{\tilde M} \subset \mc D(((A_M(\epsilon,t)+I)_0-I)^{\frac{1}{2}}).$$ It follows that
\begin{equation}
\begin{split}
   & (f, A_M(\epsilon,t) I^{\infty}(\epsilon,t)|_{\tilde M})_M\\
  = &(f, ((A_M(\epsilon, t)+I)_0-I) I^{\infty}(\epsilon, t)|_{\tilde M})_M\\ = & (f, ((A_M(\epsilon, t)+I)_0-I)^{\frac{1}{2}} ((A_M(\epsilon, t)+I)_0-I)^{\frac{1}{2}} I^{\infty}(\epsilon, t)|_{\tilde M})_M \\
  = &0.
  \end{split}
  \end{equation}
  Let $f_{\epsilon, t}$ be a function such that $f_{\epsilon, t}|_{\tilde M}=f$ and
  $$f_{\epsilon, t}(g l n)= (\mu^{\epsilon} \otimes \nu^{t+\rho})(l^{-1}) f_{\epsilon,t} (g) \qquad (l \in \tilde L, n \in N).$$
  The function $f_{\epsilon, t}$ is not necessarily in $I(\epsilon, t)$.
  By Lemma ~\ref{7.2}, $\forall \ u \in V(\epsilon, t)$,
  $$0=(f, A_M(\epsilon,t) (u|_{\tilde M}))_M=(f_{\epsilon, t}, u)_{\epsilon, t}=(f_{\epsilon, t}, A_X(\epsilon, t) u)_X.$$
  This equality is to be interpreted as an equality of integrals according to the definitions of $(\ , \ )_M$ and 
  $(\ , \ )_X$. Since $A_X(\epsilon, t)$ acts on $\tilde U(n)$-type in $V(\epsilon, t)$ as a scalar,
  $A_X(\epsilon, t)V(\epsilon, t)=V(\epsilon, t)$. We now have
  $$(f_{\epsilon, t}, V(\epsilon, t))_X=0.$$ 
  In particular, $f_{\epsilon,t}|_{\tilde U(n)} \in L^1(X)$.
Therefore $f_{\epsilon,t} =0$. We see that
$I_{A_{M}(\epsilon, t)}$ is an isometry from $C(\epsilon, t)$ onto 
$$L^2(\tilde M \times_{\widetilde{ O(q-p) U(p)} } \mathbb C_{\mu^{\epsilon}}, d [g_1] d [k_2]) .$$ 
$\Box$\\
\\
The Hilbert space $$L^2(\tilde M \times_{\widetilde{ O(q-p) U(p)} } \mathbb C_{\mu^{\epsilon}}, d [g_1] d [k_2]) $$
is the mixed model for $I(\epsilon, 0)$. We now obtain an isometry from $C(\epsilon, t)$ onto $I(\epsilon, 0)$. Denote this isometry by $\mc U(\epsilon,t)$.
Now, with in the mixed model, the action of $I(\epsilon,t)(g_1 k_2)$ is simply the left regular action and it is independent of $t$. We obtain

\begin{lem}~\label{com1} Suppose $t <0$. Let $g \in \widetilde{Sp(p, \mb R)}$ or $g \in \tilde U(q)$. Let $L(g)$ be the left regular action on
$$L^2(\tilde M \times_{\widetilde{ O(q-p) U(p)} } \mathbb C_{\mu^{\epsilon}}, d [g_1] d [k_2]) .$$ 
As a operator on $I^{\infty}(\epsilon, t)|_{\tilde M}$, $L(g)$ commutes with $A_M(\epsilon, t)$. Furthermore, $L(g)$ commutes with $(A_M(\epsilon,t)+I)_0-I$.
\end{lem}
Proof: Let $g \in \tilde M$. Both $A_M(\epsilon, t)$ and $L(g)$ are well-defined operator on $I^{\infty}(\epsilon, t)|_{\tilde M}$.
Restricting $A(\epsilon,t) I(\epsilon, t)(g)= I(\epsilon, -t)(g) A(\epsilon, t)$ onto $\tilde M$, we have
$$A_M(\epsilon, t) L(g)= L(g) A_M(\epsilon, t).$$
It follows that
$$(A_M(\epsilon, t)+I) L(g)=L(g) (A_M(\epsilon, t)+I).$$
Recall that $(A_M(\epsilon, t)+I)_0$ can be defined as the inverse of $(A_M(\epsilon, t)+I)^{-1}$, which exists and is bounded.
So $L(g)$ commutes with both $(A_M(\epsilon, t)+I)^{-1}$ and $(A_M(\epsilon, t)+I)_0$. $\Box$
\begin{lem} We have, for $g \in \tilde M$,
$\mc U(\epsilon, t) I(\epsilon,t)(g) = I(\epsilon, 0)(g) \mc U(\epsilon, t)$.
\end{lem}
Proof: It suffices to show that on the mixed model, $((\mc A_{M}(\epsilon, t)+I)_0-I)^{\frac{1}{2}}$ commutes with $L(g)$. This follows from Lemma ~\ref{com1}. $\Box$.\\
\\
Our main theorem is proved.

\end{document}